\documentclass[12pt]{article}
\usepackage{amsfonts,amsmath,amssymb}
\usepackage{color}
\def\C{\mathbb{C}}

\def\Z{\mathbb{Z}}
\def\D{\mathbb{D}}

\def\bq{ \begin{equation} }
\def\eq{ \end{equation} }
\def\ben{ \begin{eqnarray} }
\def\en{ \end{eqnarray} }

\def\frac#1#2{{#1\over #2}}

\def\on#1#2{\mathop{\vbox{\ialign{##\crcr\noalign{\kern2pt}
$\scriptstyle{#2}$\crcr\noalign{\kern2pt\nointerlineskip}
\kern-2pt$\hfil\displaystyle{#1}\hfil$\crcr}}}\limits}

\begin{document}

\baselineskip=15pt
\vspace{1cm} \centerline{\LARGE  \textbf{Deformations of complex structures on Riemann surfaces 
}}
\vskip0.4cm
\centerline{\LARGE  \textbf{and integrable structures of  Whitham type hierarchies }}

\vskip1cm \hfill
\begin{minipage}{13.5cm}
\baselineskip=10pt
{\large \bf
  A. Odesskii ${}^{}$} \\ [2ex]
{\footnotesize
${}^{}$ Brock University, St. Catharines, Canada }\\

\vskip1cm{\bf Abstract}
\bigskip

We obtain variational formulas for holomorphic objects on Riemann surfaces with respect to arbitrary local coordinates on the moduli space of complex structures. These formulas are written in terms of a canonical object on the  moduli space 
which corresponds to the pairing between the space of quadratic differentials and the tangent space to the  moduli space. This canonical object satisfies certain commutation relations which appear to be the same as the ones that  emerged in  the 
integrability theory of Whitham type hierarchies. Driven by this observation, we  develop the theory of Whitham type hierarchies integrable by hydrodynamic reductions as a theory of certain differential-geometric objects.  As an application 
we prove that the universal Whitham hierarchy is integrable by hydrodynamic reductions.

\end{minipage}

\vglue1cm \textbf{Address}:
Brock University,  Niagara Region,  500 Glenridge Ave., St. Catharines, Ont., L2S 3A1 Canada

\textbf{E-mail}:
aodesski@brocku.ca 

\newpage

\tableofcontents

\newpage

\section{Introduction}

Various calculations with holomorphic objects on a Riemann surface can be done efficiently by using    the Fay identity \cite{fey,ell}. It seems that the Fay identity contains all information about identities between the  Riemann theta function, normalized holomorphic differentials,  the 
prime form and their derivatives with respect to coordinates on a Riemann surface and on its Jacobian. On the other hand, these holomorphic objects also depend on moduli of complex structures and one needs to be able to compute variations with respect to these moduli. 
Such formulas were obtained by Rauch \cite{rau}. He represented a Riemann surface as a ramified covering of $\C P^1$  and computed variations of holomorphic objects in terms of branch points of this covering. Rauch formulas have proven its usefulness and efficiency in various 
contexts \cite{fey2,kok,shr}. It is desirable however to have universal variational formulas which are independent of a particular representation of a Riemann surface and work  for arbitrary coordinates on the  moduli space. This problem can be approached as follows. It is known that the space of  
quadratic holomorphic differentials on a Riemann surface $\mathcal{E}$ is dual to the tangent space of the moduli space $M_g$ of complex structures at the point corresponding to $\mathcal{E}$  (see \cite{kod} for a  general theory of deformations of complex structures and \cite{tate} for the  Serre duality theorem). Let $v_1,...,v_{3g-3}$ be local coordinates on $M_g$,  let
$\frac{\partial}{\partial v_1},...,\frac{\partial}{\partial v_{3g-3}}$ be the corresponding basis in the tangent space and  $g_1(p)dp^2,...,g_{3g-3}(p)dp^2$ be the dual basis in the space of quadratic differentials.  The object  $$G(p)dp^2=\sum_{i=1}^{3g-3}g_i(p)dp^2\frac{\partial}{\partial v_i}$$   does not depend on any choice of coordinates. Moreover, let $M_{g,n}$ be the moduli space of Riemann surfaces with $n$ punctures $u_1,...,u_n$. Here we can vary both the  complex structure of  $\mathcal{E}$ and points 
$u_1,...,u_n$ in $\mathcal{E}$. A basis in the tangent space looks like $\frac{\partial}{\partial u_1},...,\frac{\partial}{\partial u_n},~\frac{\partial}{\partial v_1},...,\frac{\partial}{\partial v_{3g-3}}$ and the corresponding object is
$$\hat{G}(p)dp^2=\sum_{i=1}^nF(p,u_i)dp^2\frac{\partial}{\partial u_i}+\sum_{j=1}^{3g-3}g_j(p)dp^2\frac{\partial}{\partial v_j}$$ where $F(p,u)$ has a pole of order one at $p=u$ and is  holomorphic outside the diagonal. The residue of $F(p,u)$ at $p=u$ is a 
constant and without loss of generality we assume  $F(p,u)=\frac{1}{p-u}+O(1)$. Roughly speaking, $F(p,u)$ should be a quadratic differential with respect to $p$ and a vector field with respect to $u$. However,  the  transformation law for this 
object is more complicated. Indeed, if we change coordinates by $p=\mu(\tilde{p},v_1,...,v_{3g-3})$, $u_i=\mu(\tilde{u}_i,v_1,...,v_{3g-3})$ (and do not change $v_1,...,v_{3g-3}$), then $\hat{G}(p)dp^2$ should transform as a vector field with respect to 
$u_1,...,u_n,v_1,...,v_{3g-3}$ (and as a quadratic differential with respect to $p$). This leads to the following transformation law for $F(p,u)$: 
$$\tilde{F}(\tilde{p},\tilde{u})=\frac{\mu^{\prime}(\tilde{p})^2}{\mu^{\prime}(\tilde{u})}\Big(F(\mu(\tilde{p}),\mu(\tilde{u}))-G(\mu(\tilde{p}))(\mu(\tilde{u}))\Big)$$
with the same coefficient at $F$ as if  it were a quadratic differential with respect to the first argument and a vector field with respect to the second one but with an additional term depending on $G$. 

\textit{The first main result of the paper is the set  of  explicit  formulas for the  action of $G(p)$ on various holomorphic objects on the Riemann surface  $\mathcal{E}$ such as prime form, holomorphic 1-forms, period matrix, see (\ref{GTc4}),  (\ref{GTc5}),  (\ref{GTc6}).}

We have  also obtained commutation relations for $G(p)$:

$[G(p_1),G(p_2)]=$
\begin{equation}\label{GTi}
=F(p_2,p_1)G^{\prime}(p_1)-F(p_1,p_2)G^{\prime}(p_2)+2F(p_2,p_1)_{p_1}G(p_1)-2F(p_1,p_2)_{p_2}G(p_2).
\end{equation}
The same relations hold for $\hat{G}(p)$. One can check that these relations are invariant with respect to an arbitrary change of coordinates by virtue of transformation laws of  $G$ and $F$. 

Our  main motivation for these  studies   came from attempts to understand better  integrable structures of  the so-called Whitham type hierarchies \cite{kr2,odes,odes2}. Recall that a Whitham type hierarchy is defined as compatibility conditions of the following system of PDEs:
\begin{equation}\label{whi}
\frac{\partial\psi}{\partial t_i}=h_i(z,u_1,...,u_n),~i=1,...,N. 
\end{equation}
Here $\psi,u_1,...,u_n$ are functions of times $t_1,...,t_N$ and $z$ is a parameter. The system (\ref{whi}) is understood as a parametric way of defining $N-1$ relations between 
partial derivatives $\frac{\partial\psi}{\partial t_i},~i=1,...,N$ obtained by excluding $z$ from these equations. Functions $h_i$ are called potentials of this Whitham type hierarchy.

 An important class 
of such hierarchies associated with the moduli space of Riemann surfaces of genus $g$ with $n$ punctures (the so-called universal Whitham hierarchy) was constructed 
and studied in \cite{kr2,kr1}. The universal Whitham hierarchy is important in the theory of Frobenius manifolds \cite{dub}, matrix models and other areas of mathematics. 
Note that the set of times in 
the universal Whitham hierarchy coincides with a set of meromorphic differentials on a Riemann surface (holomorphic outside punctures), and that the potentials $h_i(z)$ are 
integrals of these differentials.

A natural question is in which sense a Whitham type hierarchy is integrable. In this paper we concentrate on an approach to integrability theory of such systems based on the so-called hydrodynamic reductions \cite{fer1,fer2,fer3}.
In this approach a quasi-linear system is called integrable if it possesses a large family of  hydrodynamic reductions. This family of hydrodynamic reductions must be parametrized by solutions of another system of PDEs   called 
Gibbons-Tsarev system \cite{Gibt,os1,os11}. Therefore, Gibbons-Tsarev systems play  a  crucial role in this integrability theory. See \cite{os11} and references therein for the  definition and examples of Gibbons-Tsarev systems.

\textit{The second main result of the paper is a description of integrable structures that  appeared in the  hydrodynamic reduction method as a certain differential-geometric structure.  We  call it a GT structure.}

By definition, a GT structure is defined locally by a family of vector fields $g(p)$ and a function $f(p_1,p_2)$ satisfying  relations similar to (\ref{GTi}), see Section 3 for precise definitions.  It becomes transparent from these definitions  that a natural 
GT structure exists on the  moduli space $M_{g,n}$ and is represented by the objects $\hat{G}(p)$ and $F(p,u)$ described above.

Given a GT structure one could ask how to find all corresponding integrable Whitham type hierarchies. It turns out that in order to classify all possible integrable hierarchies with given GT structure one needs to find all functions $\lambda(p_1,p_2)$ 
satisfying  the  functional equation
$$g(p_1)(\lambda(p_2,p_3))=\lambda(p_1,p_3)\lambda(p_2,p_1)_{p_1}-\lambda(p_2,p_3)f(p_1,p_2)_{p_2}-$$
$$-f(p_1,p_2)\lambda(p_2,p_3)_{p_2}-f(p_1,p_3)\lambda(p_2,p_3)_{p_3}.$$ 
Moreover, to find all potentials $h(p)$  of a given hierarchy one needs to solve another functional equation
$$g(p_1)(h(p_2))=\lambda(p_1,p_2)h^{\prime}(p_1)-f(p_1,p_2)h^{\prime}(p_2).$$

It is natural to ask if the universal Whitham hierarchy is integrable by hydrodynamic reductions. One could expect that the corresponding GT structure is given by $\hat{G}(p)$ and $F(p_1,p_2)$ and needs to find a function $\lambda(p_1,p_2)$ which 
gives the universal Whitham hierarchy.

\textit{The third main result of this paper is a proof  that the universal Whitham hierarchy is indeed integrable in  all genera by hydrodynamic reductions.}   We give the corresponding function $\lambda(p_1,p_2)$ and  the  precise form of potentials.

Let us describe the content of the paper. In Section 2   we recall main definitions and notations of holomorphic objects on a Riemann surface, construct our main object $G(p)$ and compute its action on holomorphic objects. 
We compute commutation relations for $G(p)$ as well. We give some examples and explain how the Rauch formulas are connected with ours. In Section 3 we introduce GT structures and develop a theory of these structures. In particular, we explain 
how to construct new GT structures from a given one and how to construct potentials if we are given a function $\lambda(p_1,p_2)$ defining our hierarchy. We also explain a relation between GT structures and Lie algebroids of  a  certain type.
In Section 4 we recall  the   definition and basic properties  of  Whitham type hierarchies. In Section 5 we recall definition of Gibbon-Tsarev systems and prove that there exists a one-to-one correspondence between 
Gibbons-Tsarev systems and GT structures.
 In Section 6 we discuss  the  definition of integrability of Whitham type hierarchies in our framework.  We explain that the definition based on hydrodynamic reductions and Gibbons-Tsarev systems is equivalent to ours.
 We  refer to \cite{os11} for a  full discussion of integrability of Whitham type hierarchies based on 
hydrodynamic reductions and Gibbons-Tsarev systems. It was not possible to make this paper self-contained and repeat this discussion here without essential increasing of the  length  of the present paper.  In Section 7 we recall  the  definition of the 
universal Whitham hierarchy and prove that this hierarchy is integrable by hydrodynamic reductions.

\section{Holomorphic objects on Riemann surfaces and deformations  of complex structures}

Let $\mathcal{E}=\D/\Gamma$ be a compact Riemann surface of genus $g>1$, $\D\subset \C$ its universal covering and $\Gamma=\pi_1(\mathcal{E})$. 
Denote by $a_{\alpha},b_{\alpha},~\alpha=1,...,g$ a canonical basis in the homology group $H_1(\mathcal{E},\Z)$. Let us choose a coordinate in $\D$ and use the same symbols for 
holomorphic objects on $\mathcal{E}$ and their lifting on $\D$. We will also use the same symbol  for a point in  $\mathcal{E}$, its lifting in $\D$ and its coordinate.  Let $\omega_{\alpha}(z)dz$ be the basis of holomorphic 1-forms on 
$\mathcal{E}$ normalized by $\int_{a_{\alpha}}\omega_{\beta}dz=\delta_{\alpha\beta}$. Choose a basepoint $z_0$ and define the Abel map $q_{\alpha}(z)=\int_{z_0}^z\omega_{\alpha}(z)dz$. 
Note that $\omega_{\alpha}=q^{\prime}_{\alpha}$. Denote the prime form\footnote{In this paper we represent differential-geometric objects as functions with prescribed transformation laws  
with respect to an arbitrary change of coordinates. For example if $x=\mu(\tilde{x}),~y=\mu(\tilde{y})$, then the prime form transforms as 
$\tilde{E}(\tilde{x},\tilde{y})=\mu^{\prime}(\tilde{x})^{-1/2}\mu^{\prime}(\tilde{y})^{-1/2}E(\mu(\tilde{x}),\mu(\tilde{y}))$.}  by $E(x,y)(dx)^{-1/2}(dy)^{-1/2}$. Let $B_{\alpha\beta}= \int_{b_{\alpha}}\omega_{\beta}dz$ be the matrix of $b$-periods.   Details on holomorphic  objects on Riemann surfaces are given in  \cite{fey,ell,kok}.  Recall that 
\begin{equation}\label{id1}
 E(v,u)=-E(u,v),~~~~E(u,v)=u-v-\frac{1}{12}S(u)(u-v)^3+O((u-v)^4),
\end{equation}
where $S(p)$ is the Bergman projective connection on $\mathcal{E}$. Note that $E(u,v)$ is  multivalued. If $u$ or $v$ is moved by $a_{\alpha}$, it remains invariant.  If $u$ moves by $b_{\alpha}$ to $\bar{u}$ or 
$v$ moves by $b_{\alpha}$ to $\bar{v}$, then 
\begin{equation}\label{aut}
E(\bar{u},v)=E(u,v)\exp\Big(-\pi i B_{\alpha\alpha}+2\pi i(q_{\alpha}(v)-q_{\alpha}(u))\Big),
\end{equation}
$$E(u,\bar{v})=E(u,v)\exp\Big(-\pi i B_{\alpha\alpha}-2\pi i(q_{\alpha}(v)-q_{\alpha}(u))\Big).$$

Let $W(u,v)=(\ln(E(u,v))_{uv}$ be the Bergman kernel. Recall that 
\begin{equation}\label{aut1}
\int_{a_i}W(u,v) du=0,~\int_{b_{\alpha}}W(u,v) du=2\pi i\omega_{\alpha}(v),~ \int_{b_{\alpha}}\int_{b_{\beta}}W(u,v) du dv=2\pi iB_{\alpha\beta}.
\end{equation}

Recall a description of  the  tangent space to the moduli space $M_g$ of Riemann surfaces at the point corresponding to  $\mathcal{E}$ \cite{kont,beil}.
Let $p\in\mathcal{E}$  be the center of  a small disc  $D\subset\mathcal{E}$.  Let $L$ be the Lie algebra of holomorphic vector fields on $D \setminus  \{p\}$ and $L_p$, $L_{out}$ be  subalgebras of $L$ consisting of vector fields holomorphic 
at $p$ and holomorphic on $\mathcal{E}\setminus\{p\}$ correspondingly. It is known that the tangent space to the moduli space $M_g$  is isomorphic to the quotient $L/(L_p\oplus L_{out})$. Let $M_{g,1}$ be the moduli space of Riemann surfaces with 
a puncture at $u\in\mathcal{E}$. The tangent space to $M_{g,1}$ is isomorphic to the quotient $L/(L_p\oplus L_{out,u})$ where $L_{out,u}\subset L_{out}$ consists of vector fields with zero at $u$. Let us construct vector spaces dual to these tangent  
spaces using the Serre duality theorem \cite{tate}. There exists a non degenerate pairing between the space $L$ and the space $Q$ of quadratic differentials holomorphic on  $D \setminus  \{p\}$. This pairing is given by $(v,q)=\text{Res}_p(vq)$. The space dual to the 
tangent space of $M_g$ is equal to $(L_p\oplus L_{out})^{\perp}\subset Q$ and consists of quadratic differentials holomorphic on $\mathcal{E}$. Similarly, the space dual to the 
tangent space of $M_{g,1}$ is equal to $(L_p\oplus L_{out,u})^{\perp}\subset Q$ and consists of quadratic differentials holomorphic on $\mathcal{E}\setminus\{u\}$ with pole of order less or equal to  one  at $u$. More generally, the space dual to the 
tangent space of $M_{g,n}$ of the moduli space of Riemann surfaces with punctures at $u_1,...,u_n$ consists of quadratic differentials holomorphic on $\mathcal{E}\setminus\{u_1,...,u_n\}$ with poles of order less or equal  to  one  at $u_1,...,u_n$.

Let $v_1,...,v_{3g-3}$ be local coordinates on moduli space $M_g$. Let $\frac{\partial}{\partial v_1},...,\frac{\partial}{\partial v_{3g-3}}$ be the corresponding basis in the tangent space and 
$g_1(p)dp^2,...,g_{3g-3}(p)dp^2$ be  the dual basis in the space of quadratic differentials. The object\footnote{Note that the functions $G$, $g_i$, $F$ etc. depend also on $v_1,...,v_{3g-3}$. We will often omit these arguments in order to simplify formulas.}
 $$G(p)dp^2=\sum_{i=1}^{3g-3}g_i(p)dp^2\frac{\partial}{\partial v_i}$$  does not depend on the choice  
of coordinates. A  similar construction for $M_{g,n}$ gives the object $$\hat{G}(p)dp^2=\sum_{i=1}^nF(p,u_i)dp^2\frac{\partial}{\partial u_i}+\sum_{j=1}^{3g-3}g_j(p)dp^2\frac{\partial}{\partial v_j}$$  where $u_1,...,u_n$ are coordinates of 
$n$ points in $\mathcal{E}$ and 
\begin{equation}\label{diag}
F(p_1,p_2)=\frac{1}{p_1-p_2}+O(1).
\end{equation}

{\bf Proposition  2.1.} Under an arbitrary change of coordinates of the form  
\begin{equation}\label{CHcoor}
p=\mu(\tilde{p},v_1,...,v_{3g-3}), ~u_i=\mu(\tilde{u}_i,v_1,...,v_{3g-3})
\end{equation}
the objects $G(p)$, $F(p_1,p_2)$ obey the following transformation rules

\begin{equation}\label{CHc}
\tilde{G}(\tilde{p})=\mu^{\prime}(\tilde{p})^2G(\mu(\tilde{p})),
\end{equation}
\begin{equation}\label{CHc1}
\tilde{F}(\tilde{p}_1,\tilde{p}_2)=\frac{\mu^{\prime}(\tilde{p}_1)^2}{\mu^{\prime}(\tilde{p}_2)}\Big(F(\mu(\tilde{p}_1),\mu(\tilde{p}_2))-G(\mu(\tilde{p}_1))(\mu(\tilde{p}_2))\Big)
\end{equation}

{\bf Proof.} The relation (\ref{CHc}) means that $G(p)$ is a quadratic differential in $p$ (with values in vector fields in $v_1,...,v_{3g-3}$). In order to obtain  (\ref{CHc1}) we  perform  an arbitrary change of 
coordinates of the form $p=\mu(\tilde{p},v_1,...,v_{3g-3}),~u_i=\mu(\tilde{u}_i,v_1,...,v_{3g-3}),~v_j=\tilde{v}_j$ and require that the object  $\hat{G}(p)$  transforms as a vector field in $u_1,...,u_n,v_1,...,v_{3g-3}$. The relation $(\ref{CHc1})$  is a consequence of this  requirement.   $\Box$

{\bf Proposition 2.2.} The following identities hold

\begin{equation}\label{GTc1}
[G(p_1),G(p_2)]=F(p_2,p_1)G^{\prime}(p_1)-F(p_1,p_2)G^{\prime}(p_2)+
\end{equation}
$$+2F(p_2,p_1)_{p_1}G(p_1)-2F(p_1,p_2)_{p_2}G(p_2),$$

\begin{equation}\label{GTc2}
[\hat{G}(p_1),\hat{G}(p_2)]=F(p_2,p_1)\hat{G}^{\prime}(p_1)-F(p_1,p_2)\hat{G}^{\prime}(p_2)+
\end{equation}
$$+2F(p_2,p_1)_{p_1}\hat{G}(p_1)-2F(p_1,p_2)_{p_2}\hat{G}(p_2),$$

\begin{equation}\label{GTc3}
G(p_2)(F(p_1,p_3))-G(p_1)(F(p_2,p_3))=F(p_1,p_2)F(p_2,p_3)_{p_2}-F(p_2,p_1)F(p_1,p_3)_{p_1}+
\end{equation}
$$+F(p_1,p_3)F(p_2,p_3)_{p_3}-F(p_2,p_3)F(p_1,p_3)_{p_3}+2F(p_2,p_3)F(p_1,p_2)_{p_2}-2F(p_1,p_3)F(p_2,p_1)_{p_1},$$

\begin{equation}\label{GTc4}
\frac{G(p_1)(E(p_2,p_3))}{E(p_2,p_3)}=\frac{1}{2}F(p_1,p_2)_{p_2}+\frac{1}{2}F(p_1,p_3)_{p_3}-
\end{equation}
$$-F(p_1,p_2)\frac{E(p_2,p_3)_{p_2}}{E(p_2,p_3)}-F(p_1,p_3)\frac{E(p_2,p_3)_{p_3}}{E(p_2,p_3)}-\frac{1}{2}\Big(\frac{E(p_1,p_2)_{p_1}}{E(p_1,p_2)}-\frac{E(p_1,p_3)_{p_1}}{E(p_1,p_3)}\Big)^2,$$

\begin{equation}\label{GTc5}
G(p_1)\Big(\int_{p_2}^{p_3}\omega_i\Big)=F(p_1,p_2)\omega_i(p_2)-F(p_1,p_3)\omega_i(p_3)-
\end{equation}
$$-\frac{E(p_1,p_2)_{p_1}}{E(p_1,p_2)}\omega_i(p_1)+\frac{E(p_1,p_3)_{p_1}}{E(p_1,p_3)}\omega_i(p_1),$$

\begin{equation}\label{GTc6}
G(p)(B_{jk})=2\pi i\omega_j(p)\omega_k(p).
\end{equation}

{\bf Proof.} Notice that  (\ref{GTc2}) is a formal consequence of  (\ref{GTc1}) and  (\ref{GTc3}) (see Proposition 3.1).

Consider the difference of the l.h.s. and the r.h.s. of each of  (\ref{GTc1}),  (\ref{GTc3}),  (\ref{GTc4}),  (\ref{GTc5}). Expanding these expressions on each diagonal $p_i=p_j,~i\ne j$ and using (\ref{id1}) and (\ref{diag}) one can check 
that each of these expressions is  holomorphic on all diagonals. Making an arbitrary change of coordinates of the form $p_i=\mu(\tilde{p_i},v_1,...,v_{3g-3}), ~i=1,2,3$ one can check that all these differences are transformed as tensor fields in 
$p_1,p_2,p_3$.   In particular, the difference between the l.h.s. and the r.h.s. of (\ref{GTc3}) is a holomorphic quadratic differential in $p_1,p_2$ and holomorphic vector field in $p_3$. This proves (\ref{GTc3}) because any holomorphic vector field vanishes. 
Similarly, the differences between the l.h.s. and the r.h.s. of (\ref{GTc4}), (\ref{GTc5}) are holomorphic quadratic differentials in $p_1$ and holomorphic functions in $p_2,p_3$. Moreover, these functions vanish on the diagonal $p_2=p_3$. This would 
prove (\ref{GTc4}), (\ref{GTc5})  (any holomorphic function is a constant) provided that we prove that the differences between the l.h.s. and the r.h.s.  are single valued. 

Taking the second derivative of the equation (\ref{GTc4}) we get
\begin{equation}\label{GTc41}
G(p_1)(W(p_2,p_3))=
\end{equation}
$$=-\Big(F(p_1,p_2)\frac{E(p_2,p_3)_{p_2}}{E(p_2,p_3)}+F(p_1,p_3)\frac{E(p_2,p_3)_{p_3}}{E(p_2,p_3)}-\frac{E(p_1,p_2)_{p_1}E(p_1,p_3)_{p_1}}{E(p_1,p_2)E(p_1,p_3)}\Big)_{p_2p_3}$$
where $(W(p_2,p_3)=(\ln(E(p_2,p_3))_{p_2p_3}$ is the Bergman kernel. Let  us  prove this identity.  Let $\Delta(p_1,p_2,p_3)$ be the difference of the l.h.s. and the r.h.s. of  (\ref{GTc41}). It is a quadratic differential in 
$p_1$ and 1-form in both $p_2,p_3$.  Using transformation 
properties  (\ref{aut}) we see that $\Delta(p_1,p_2,p_3)$ is single valued. Therefore, $\Delta(p_1,p_2,p_3)=\sum_{\alpha,\beta=1}^gr_{\alpha\beta}(p_1)\omega_{\alpha}(p_2)\omega_{\beta}(p_3)$ where $r_{\alpha\beta}(p_1)$ are 
some holomorphic quadratic differentials. Computing $\int_{a_{\alpha}}\int_{a_{\beta}}\Delta(p_1,p_2,p_3)dp_2 dp_3$ we obtain $r_{\alpha\beta}(p_1)=0$ which proves  (\ref{GTc41}). Computing $\int_{b_{\alpha}}\int_{b_{\beta}}~dp_2 dp_3$ 
of the l.h.s. and the r.h.s. of  (\ref{GTc41}) and using (\ref{aut1}) and  (\ref{aut}) we obtain (\ref{GTc6}). The difference between the l.h.s.  and the r.h.s.  of (\ref{GTc5}) is single valued by virtue of  (\ref{GTc6}). This proves (\ref{GTc5}). 
Equation (\ref{GTc4}) is proven in 
 a  similar way. Note that the difference between the l.h.s. and the r.h.s. of (\ref{GTc4}) is single valued by virtue of (\ref{GTc5}). Equation (\ref{GTc1}) is proven by applying its  l.h.s.  and the  r.h.s. to $B_{jk}$. For example,  on the l.h.s. we have 
  $G(p_1)(G(p_2)(B_{jk}))-G(p_2)(G(p_1)(B_{jk}))$.
Computing by virtue of  (\ref{GTc6}),  (\ref{GTc5}) we prove  (\ref{GTc1}). $\Box$

{\bf Remark 2.1.} Recall that the Riemann theta-function is defined by
$$\theta(z_1,...,z_g)=\sum_{{\bf m}\in\Z^g}\exp(2\pi i{\bf m}\cdot{\bf z}+\pi i{\bf m}{\bf B}{\bf m}^t).$$ 
 Here  we use bold  symbols for the corresponding vectors: ${\bf m}=(m_1,...,m_g),~{\bf z}=(z_1,...,z_g)$, ${\bf m}\cdot{\bf z}=m_1z_1+...+m_gz_g$, and ${\bf B}$ is the period matrix. We have 
$$G(p)(\theta(z_1,...,z_g))=\sum_{\alpha,\beta=1}^g\frac{\partial\theta(z_1,...,z_g)}{\partial B_{\alpha\beta}}G(p)(B_{\alpha\beta})=
\frac{1}{2}\sum_{\alpha,\beta=1}^g\frac{\partial^2\theta(z_1,...,z_g)}{\partial z_{\alpha}\partial z_{\beta}}\omega_{\alpha}(p)\omega_{\beta}(p)$$
where we used heat equation for $\theta$ and  (\ref{GTc6}).

{\bf Remark 2.2.} Expanding  (\ref{GTc4}) on diagonal $p_2=p_3$ we obtain
$$G(p_1)(S(p_2))+F(p_1,p_2)_{p_2^3}+2S(p_2)F(p_1,p_2)_{p_2}+S(p_2)_{p_2}F(p_1,p_2)-6W(p_1,p_2)^2=0.$$

{\bf Example 2.1. } Let $g=2$. Represent $\mathcal{E}$ as a 2-fold covering of $\C P^1$. Let $x$ be an affine coordinate in $\C P^1$ and  let  branch points of the covering  be  at $x=0,1,\infty,a,b,c$. The curve $\mathcal{E}$ is given by 
$y^2=x(x-1)(x-a)(x-b)(x-c)$. One can check that 
$$G(p)=\frac{1}{2p(p-1)}\Big(\frac{a(a-1)}{p-a}\frac{\partial}{\partial a}+\frac{b(b-1)}{p-b}\frac{\partial}{\partial b}+\frac{c(c-1)}{p-c}\frac{\partial}{\partial c}\Big),$$
$$F(p_1,p_2)=\frac{(p_1-a)(p_1-b)(p_1-c)p_2(p_2-1)+q_1q_2}{2(p_1-p_2)p_1(p_1-1)(p_1-a)(p_1-b)(p_1-c)}$$
where $p,p_1,p_2$ are affine coordinates in $\C P^1$ and $q_i^2=p_i(p_i-1)(p_i-a)(p_i-b)(p_i-c)$.

Let us compare our variational formulas with Rauch ones. The equation (\ref{GTc6}) reads 
\begin{equation}\label{rau}
\frac{1}{2p(p-1)}\Big(\frac{a(a-1)}{p-a}\frac{\partial}{\partial a}+\frac{b(b-1)}{p-b}\frac{\partial}{\partial b}+\frac{c(c-1)}{p-c}\frac{\partial}{\partial c}\Big)(B_{jk})dp^2=2\pi i\omega_j(p)dp\cdot\omega_k(p)dp.
\end{equation}
Let $\tau$ be a local coordinate near branch point $a$, we have $p=a+\tau^2$, $dp=2\tau d\tau$. Expanding the l.h.s. of  (\ref{rau}) we get  $(2\frac{\partial B_{jk}}{\partial a}+O(\tau))d\tau^2$. Therefore 
$\frac{\partial B_{jk}}{\partial a}=\pi i \frac{\omega_j(p)dp}{d\tau}|_{p=a}\frac{\omega_k(p)dp}{d\tau}|_{p=a}$ and we arrive  at  a Rauch formula. 

In general, if $\mathcal{E}$ is represented as a branched covering of $\C P^1$ ramified at $a_k\in \C P^1$ with ramification indexes $r_k$, $k=1,2,...$, then $G(p)dp^2=\Big(\frac{1}{r_k(p-a_k)}\frac{\partial}{\partial a_k}+o((p-a_k)^{-1})\Big)dp^2$ and Rauch 
formulas can be derived from ours in a similar way.

{\bf Example 2.2.} Let us choose $B_{j_1,k_1},...,B_{j_{3g-3},k_{3g-3}}$ as local coordinates  in $M_g$. Applying (\ref{GTc6}) to $B_{j_l,k_l}$ we get $$G(p)=2\pi i \sum_{l=1}^{3g-3}\omega_{j_l}(p)\omega_{k_l}(p)\frac{\partial}{\partial B_{j_l,k_l}}.$$
Applying again  (\ref{GTc6}) to an arbitrary $B_{jk}$ we obtain quadratic relations between normalized differentials. Namely, if $S(B_{11},B_{12},...,B_{gg})=0$ is a relation between the entries of the period matrix (there are $\frac{g(g+1)}{2}-3g+3$ 
functionally independent ones), then $$\sum_{j,k=1}^g\frac{\partial S}{\partial B_{jk}}\omega_j(p)\omega_k(p)=0.$$  See \cite{andr} for a  similar  formula  for quadratic relations between normalized differentials.

\section{GT structures}

Based on the identities  (\ref{GTc1}),  (\ref{GTc3}) we want to introduce a general differential-geometric structure: a family of vector fields $g(p)$ and a function $f(p_1,p_2)$ satisfying the same relations. We will see later that this structure is 
equivalent to an integrability structure of  Whitham type hierarchies, the so-called Gibbons-Tsarev system.

Let $g(p)=\sum_{i=1}^{m}g_i(p,v_1,...v_m)\frac{\partial}{\partial v_i}$ be a family of vector fields parameterized by $p$ and $f(p_1,p_2,v_1,...,v_m)$ be a function.

{\bf Definition 3.1.} A local GT structure is a family $g(p)$ and a function $f(p_1,p_2)$ satisfying the following relations:

\begin{equation}\label{GT1}
[g(p_1),g(p_2)]=f(p_2,p_1)g^{\prime}(p_1)-f(p_1,p_2)g^{\prime}(p_2)+2f(p_2,p_1)_{p_1}g(p_1)-2f(p_1,p_2)_{p_2}g(p_2),
\end{equation}

\begin{equation}\label{GT2}
g(p_2)(f(p_1,p_3))-g(p_1)(f(p_2,p_3))=f(p_1,p_2)f(p_2,p_3)_{p_2}-f(p_2,p_1)f(p_1,p_3)_{p_1}+
\end{equation}
$$+f(p_1,p_3)f(p_2,p_3)_{p_3}-f(p_2,p_3)f(p_1,p_3)_{p_3}+2f(p_2,p_3)f(p_1,p_2)_{p_2}-2f(p_1,p_3)f(p_2,p_1)_{p_1},$$

\begin{equation}\label{GT3}
f(p_1,p_2)=\frac{1}{p_1-p_2}+O(1).
\end{equation}

Here and in the sequel we often omit additional arguments $v_1,...v_m$, indexes stand for partial derivatives and $g^{\prime}(p)=\frac{\partial g(p,v_1,...,v_m)}{\partial p}$.

Given a GT structure we can construct new GT structures in different ways. 

{\bf Proposition 3.1.} Let $g(p)$, $f(p_1,p_2)$ satisfy relations (\ref{GT1}),  (\ref{GT2}) and 

\begin{equation}\label{ext}
\hat{g}(p)=f(p,u_1)\frac{\partial}{\partial u_1}+...+f(p,u_n)\frac{\partial}{\partial u_n}+g(p).
\end{equation}

Then $\hat{g}(p)$, $f(p_1,p_2)$ also satisfy  relations (\ref{GT1}),  (\ref{GT2}).

{\bf Proof.} Equation  (\ref{GT1}) is verified by direct computation for $n=1$ and through  induction by $n$ for $n>1$. Equation  (\ref{GT2}) remains the same because $f(p_1,p_2)$ does not depend on $u_1,...,u_n$. $\Box$

We say that  a  GT structure given by $\hat{g}(p)$, $f(p_1,p_2)$ is obtained from a  GT structure $g(p)$, $f(p_1,p_2)$ by adding $n$ points $u_1,...,u_n$. This procedure corresponds to a  regular fields extension of a Gibbons-Tsarev system \cite{os11}.   

{\bf Proposition 3.2.} Let $g(p)$, $f(p_1,p_2)$ satisfy relations (\ref{GT1}),  (\ref{GT2}) and 

\begin{equation}\label{ext1}
\hat{g}^{(n_1,...,n_k)}(p)=
\end{equation}
$$=\sum_{\substack{1\leq j\leq k, \\ 0\leq i_{j,1},...,i_{j,n_j}, \\ i_{j,1}+...+i_{j,n_j}\leq n_{j}}}\frac{(i_{j,1}+2i_{j,2}+...+n_ji_{j,n_j})!}{i_{j,1}!...i_{j,n_j}!~1!^{i_{j,1}}...n_j!^{i_{j,n_j}}}\frac{\partial^{i_{j,1}+...+i_{j,n_j}}f(p,u_{j,0})}
{\partial u_{j,0}^{i_{j,1}+...+i_{j,n_j}}}\frac{\partial}{\partial u_{j,i_{j,1}+...+i_{j,n_j}}}+g(p).$$

Then $\hat{g}^{(n_1,...,n_k)}(p)$, $f(p_1,p_2)$ also satisfy  relations (\ref{GT1}),  (\ref{GT2}).

{\bf Proof.}  Let us start with the following local GT structure

\begin{equation}\label{ext2}
\hat{g}(p)=\sum_{\substack{1\leq j\leq k, \\ 0\leq l\leq n_j}}f(p,v_{j,l})\frac{\partial}{\partial v_{j,l}}+g(p).
\end{equation}

We make the following change of coordinates

$$v_{j,0}=u_{j,0},$$
$$v_{j,1}=u_{j,0}+\epsilon u_{j,1},$$
\begin{equation}\label{col}
v_{j,2}=u_{j,0}+2\epsilon u_{j,1}+\epsilon^2 u_{j,2},
\end{equation}
$$.............$$
$$v_{j,n_j}=u_{j,0}+n_j\epsilon u_{j,1}+\frac{n_j(n_j-1)}{2}\epsilon^2u_{j,2}+...+\epsilon^{n_j}u_{j,n_j}.$$

In new coordinates we have
$$\hat{g}(p)=\sum_{1\leq j\leq k}\Big(f(p,u_{j,0})\frac{\partial}{\partial u_{j,0}}+\frac{1}{\epsilon}(f(p,u_{j,0}+\epsilon u_{j,1})-f(p,u_{j,0}))\frac{\partial}{\partial u_{j,1}}+$$
$$+\frac{1}{\epsilon^2}(f(p,u_{j,0}+2\epsilon u_{u,1}+\epsilon^2u_{j,2})-2f(p,u_{j,0}+\epsilon u_{j,1})+f(p,u_{j,0}))\frac{\partial}{\partial u_{j,2}}+...$$
$$+\frac{1}{\epsilon^{n_j}}(f(p,u_{j,0}+n_j\epsilon u_{j,1}+...+\epsilon^{n_j}v_{j,n_j})-...+(-1)^{n_j}f(p,u_{j,0}))\frac{\partial}{\partial u_{j,n_j}})\Big)+g(p).$$

Taking the limit $\epsilon\to 0$ we obtain  (\ref{ext1}). $\Box$

We say that the GT structure  (\ref{ext1}) is obtained from the GT structure  (\ref{ext2}) by colliding  points  $v_{j,0},v_{j,1},...,v_{j,n_j}$ for each $j$.

{\bf Remark 3.1.} Equation  (\ref{GT2}) is equivalent to Jacobi identity for  (\ref{GT1}) provided that vector fields $g(p_1),~g(p_2),~g(p_3),~g^{\prime}(p_1),~g^{\prime}(p_2),~g^{\prime}(p_3)$ are linearly  
independent for generic $p_1,~p_2,~p_3$.

{\bf Remark 3.2.} A local GT structure can be regarded as a certain Lie algebroid. Let $$g(p)=e_2+(p-z)e_3+(p-z)^2e_4+....$$ In other words, let $e_{i+2}=i!~g^{(i)}(z)$. Let $e_1=\frac{\partial}{\partial z}$ and 
$$f(p_1,p_2)=\frac{1}{p_1-p_2}+\sum_{i,j=0}^{\infty}f_{i,j}(z)~(p_1-z)^i(p_2-z)^j.$$
Then we have $[e_1,e_i]=(i-1)e_{i+1}$ and equation  (\ref{GT1}) is equivalent to
$$[e_i,e_j]=(j-i)e_{i+j}+\sum_{r=0}^{i-1}(i+r-1)f_{j-2,r}e_{i-r+1}-\sum_{r=0}^{j-1}(j+r-1)f_{i-2,r}e_{j-r+1}.$$
In particular, if $f(p_1,p_2)=\frac{1}{p_1-p_2}$, then we get $[e_i,e_j]=(j-i)e_{i+j}$ for $e_1,e_2,...$. Note that  (\ref{GT2}) always holds for  $f(p_1,p_2)=\frac{1}{p_1-p_2}$. Therefore, a  local GT structure can be regarded as a certain  deformation of a Lie 
algebra with basis $e_1,e_2,...$ and bracket  $[e_i,e_j]=(j-i)e_{i+j}$ in  the  class of Lie algebroids.

Given a local GT structure one wants to classify all Whitham type hierarchies  that are  integrable by hydrodynamic reductions and  that  correspond to a given Gibbons-Tsarev system. It turns out that in order to do this one needs to find all functions $\lambda(p_1,p_2,v_1,...,v_m)$ 
satisfying a certain condition. This can be formalized  in the following way:

{\bf Definition 3.2.} An enhanced local GT structure is a family of vector fields $g(p)$, a function $f(p_1,p_2)$ and an additional function $\lambda(p_1,p_2,v_1,...,v_m)$ satisfying the  relations  (\ref{GT1}),  (\ref{GT2}),  (\ref{GT3})  and 

\begin{equation}\label{GT4}
g(p_1)(\lambda(p_2,p_3))=\lambda(p_1,p_3)\lambda(p_2,p_1)_{p_1}-\lambda(p_2,p_3)f(p_1,p_2)_{p_2}-
\end{equation}
$$-f(p_1,p_2)\lambda(p_2,p_3)_{p_2}-f(p_1,p_3)\lambda(p_2,p_3)_{p_3},$$

$$\lambda(p_1,p_2)=\frac{1}{p_1-p_2}+O(1).$$

Given an  enhanced local GT structure one wants to find a vector space of  all potentials of the corresponding Whitham type hierarchy. In all known examples these spaces are spaces of solutions of linear systems of PDEs. However, in the  general case 
we can define this vector space as a space of solutions of a  linear functional equation.

{\bf Definition 3.3.} Given an  enhanced local GT structure we define the corresponding vector space of potentials as the space of solutions of the following functional equation  for a function $h(p,v_1,...,v_m)$:

\begin{equation}\label{GT5}
g(p_1)(h(p_2))=\lambda(p_1,p_2)h^{\prime}(p_1)-f(p_1,p_2)h^{\prime}(p_2).
\end{equation}

Note that expanding (\ref{GT5}) near diagonal $p_2=p_1$ we obtain for $h(p,v_1,...,v_m)$  a system of linear PDEs equivalent to  (\ref{GT5}).

The following procedure gives a standard way to obtain solutions of  (\ref{GT5}):

{\bf Proposition 3.3.} Let $\gamma$ be a path in $\C$ such that $\int_{\gamma}\frac{\partial(\lambda(t,p_2)f(p_1,t))}{\partial t}d t=0$. Then $$h(p)=\int_{\gamma}\lambda(t,p)d t$$  is 
a solution of  (\ref{GT5}).

{\bf Proof.} Substitute this expression for $h(p)$ into  (\ref{GT5}) and use  (\ref{GT4}).   Direct computation shows that the difference between the r.h.s and the l.h.s. of  (\ref{GT5}) is $\int_{\gamma}\frac{\partial(\lambda(t,p_2)f(p_1,t))}{\partial t}d t$. $\Box$

Let us promote local GT structures to differential-geometric ones. 

{\bf Proposition 3.4.} Relations  (\ref{GT1}),  (\ref{GT2}),  (\ref{GT3}) are invariant with respect to arbitrary transformations of the form 

\begin{equation}\label{CH1}
p_i=\mu(\tilde{p}_i,v_1,...,v_m),~~~~~~~\tilde{g}(\tilde{p})=\mu^{\prime}(\tilde{p})^2g(\mu(\tilde{p})),
\end{equation}
$$\tilde{f}(\tilde{p}_1,\tilde{p}_2)=\frac{\mu^{\prime}(\tilde{p}_1)^2}{\mu^{\prime}(\tilde{p}_2)}\Big(f(\mu(\tilde{p}_1),\mu(\tilde{p}_2))-g(\mu(\tilde{p}_1))(\mu(\tilde{p}_2))\Big).$$

Let $\pi:~M\to B$ be a bundle with $m$ dimensional fiber $F$ and one dimensional base $B$. 

{\bf Definition 3.4.} A GT structure on $\pi$ is a local GT structure on each trivialization for  each $U\subset B$ such that for different trivializations these local GT structures are connected by  (\ref{CH1}). Here $v_1,...,v_m$ stands for coordinates on $F$ 
and $p$ is a coordinate on $B$.

{\bf Proposition 3.5.}  Relations  (\ref{GT4}) are invariant with respect to an  arbitrary transformations of the form   (\ref{CH1}) provided that $\lambda$ is transformed as

\begin{equation}\label{CH2}
\tilde{\lambda}(\tilde{p}_1,\tilde{p}_2)=\mu^{\prime}(\tilde{p}_1)\lambda(\mu(\tilde{p}_1),\mu(\tilde{p}_2))
\end{equation}

{\bf Definition 3.5.} An  enhanced  GT structure on $\pi$ is an  enhanced local GT structure on  each trivialization for each $U\subset B$ such that for different trivializations these  enhanced local GT structures are connected by  (\ref{CH1}),  (\ref{CH2}).

{\bf Example 3.1.} It is clear from (\ref{GTc1}), (\ref{GTc3})  that $g(p)=G(p),~f(p_1,p_2)=F(p_1,p_2)$ is a GT structure on the bundle $M_{g,1}\to M_g$. 

Similar GT structures exist for $g=0,1$. In the case $g=0$ we consider  the  moduli space $M_{0,n+3}$  of complex structures on $\C P^1$ with punctures in $n+3$ points. We fix 3 points at $0,1,\infty$ and move other points. The formulas for the corresponding GT structures read
\begin{equation}\label{g0}
 f(p_1,p_2)=\frac{p_2(p_2-1)}{(p_1-p_2)p_1(p_1-1)},~~~g(p)=\sum_{i=1}^{n}\frac{u_i(u_i-1)}{(p-u_i)p(p-1)}\frac{\partial}{\partial u_i}.
\end{equation}

 In the case $g=1$ we consider the  moduli space $M_{1,n+1}$   of complex structures on an elliptic curve with punctures in $n+1$ points.  We fix one  point at $0$ and move other points. We also deform  the  complex structure on our elliptic curve. The space of complex structures is one dimensional in this case. We use  the  modular parameter $\tau$ with $\text{Im}  \tau>0$ as a coordinate on the moduli space of elliptic curves. The formulas for the corresponding GT structures read
\begin{equation}\label{g1}
f(p_1,p_2)=\rho(p_1-p_2,\tau)-\rho(p_1),~~~g(p)=2\pi i\frac{\partial}{\partial\tau}+\sum_{j=1}^n(\rho(p-u_j,\tau)-\rho(p,\tau))\frac{\partial}{\partial u_j}
\end{equation}
where
$\rho(p,\tau)=\frac{\partial}{\partial p}\ln(\theta(p,\tau))$ and $\theta(p,\tau)=\sum_{k\in\Z}(-1)^ke^{2\pi i(kp+\frac{k(k-1)}{2}\tau)}.$

{\bf Remark 3.3.} In these GT structures we can also collide points and obtain new GT structures. Moreover, in the case $g=0$ (resp. $g=1$) we can collide points with $0,1,\infty$ (resp. with 0) by doing a substitution similar to (\ref{col}). In the case $g=0$ we 
can also make an arbitrary fractional linear transformation with constant coefficients sending $0,1,\infty$ to $a,b,c$ and collide some of $a,b,c$.

{\bf Remark 3.4.} Consider an enhanced local GT structure with $g(p)$ given by  (\ref{ext2}). Colliding points $v_{j,0},v_{j,1},...,v_{j,n_j}$ by substitution  (\ref{col}) and taking the limit $\epsilon\to 0$ we can do the same substitution and limit in the 
function $\lambda$ and obtain a new enhanced local GT structure. 

\section{Whitham type hierarchies}

Given a set of independent variables $t_1,...,t_N$ called times, a set of dependent variables $v_1,...,v_m$ called fields and a set of functions $h_i(z,v_1,...,v_m),~i=1,...,N$ 
called potentials we define a Whitham type hierarchy as compatibility conditions of the following system of PDEs:
\begin{equation}\label{wh}
\frac{\partial\psi}{\partial t_i}=h_i(z,v_1,...,v_m),~i=1,...,N. 
\end{equation}
Here $\psi,v_1,...,v_m$ are functions of times $t_1,...,t_N$ and $z$ is a parameter. The system (\ref{wh}) is understood as a parametric way of defining $N-1$ relations between 
partial derivatives $\frac{\partial\psi}{\partial t_i},~i=1,...,N$ obtained by eliminating  $z$ from these equations. Let us assume that the system (\ref{wh}) is compatible. Compatibility 
conditions can be written as 
\begin{equation}\label{comp}
\sum_{l=1}^m\Big(\Big(\frac{\partial h_i}{\partial z}\frac{\partial h_j}{\partial v_l}-
\frac{\partial h_j}{\partial z}\frac{\partial h_i}{\partial v_l}\Big)\frac{\partial v_l}{\partial t_k}+\Big(\frac{\partial h_j}{\partial z}\frac{\partial h_k}{\partial v_l}-
\frac{\partial h_k}{\partial z}\frac{\partial h_j}{\partial v_l}\Big)\frac{\partial v_l}{\partial t_i}+\Big(\frac{\partial h_k}{\partial z}\frac{\partial h_i}{\partial v_l}-
\frac{\partial h_i}{\partial z}\frac{\partial h_k}{\partial v_l}\Big)\frac{\partial v_l}{\partial t_j}\Big)=0 
\end{equation}
where $i,~j,~k=1,...,N$ are pairwise distinct. Let $V_{i,j,k}$ be the  linear space of functions in $z$ spanned by 
$\frac{\partial h_i}{\partial z}\frac{\partial h_j}{\partial v_l}-
\frac{\partial h_j}{\partial z}\frac{\partial h_i}{\partial v_l},~\frac{\partial h_j}{\partial z}\frac{\partial h_k}{\partial v_l}-
\frac{\partial h_k}{\partial z}\frac{\partial h_j}{\partial v_l},~\frac{\partial h_k}{\partial z}\frac{\partial h_i}{\partial v_l}-
\frac{\partial h_i}{\partial z}\frac{\partial h_k}{\partial v_l},~l=1,...,m$.

{\bf Proposition  4.1.} Let $V_{i,j,k}$ be finite dimensional and  $\dim V_{i,j,k}=D$. Then  (\ref{comp}) is equivalent to a hydrodynamic type system of  $D$  linearly independent 
equations of the form
\begin{equation}\label{hydr}
\sum_{l=1}^m\Big(a_{rl}(v_1,...,v_m)\frac{\partial v_l}{\partial t_i}+b_{rl}(v_1,...,v_m)\frac{\partial v_l}{\partial t_j}
+c_{rl}(v_1,...,v_m)\frac{\partial v_l}{\partial t_k}\Big)=0,~r=1,..., D.
\end{equation}

{\bf Proof.} Let $\{S_1(z),...,S_D(z)\}$ be a basis in $V_{i,j,k}$ and 

$\frac{\partial h_i}{\partial z}\frac{\partial h_j}{\partial v_l}-
\frac{\partial h_j}{\partial z}\frac{\partial h_i}{\partial v_l}=\sum_{r=1}^Dc_{rl}S_r,~\frac{\partial h_j}{\partial z}\frac{\partial h_k}{\partial v_l}-
\frac{\partial h_k}{\partial z}\frac{\partial h_j}{\partial v_l}=\sum_{r=1}^Da_{rl}S_r,~\frac{\partial h_k}{\partial z}\frac{\partial h_i}{\partial v_l}-
\frac{\partial h_i}{\partial z}\frac{\partial h_k}{\partial v_l}=\sum_{r=1}^Db_{rl}S_r.$ 

Substituting these expressions into (\ref{comp}) and equating to zero coefficients at $S_1,...,S_D$ we obtain (\ref{hydr}). $\Box$

{\bf Remark 4.1.} In the theory of integrable systems of hydrodynamic type the system (\ref{wh}) is often referred to as a pseudo-potential representation of the system (\ref{hydr}).

{\bf Remark 4.2.} In all known examples of integrable Whitham type hierarchies we have $n\leq D\leq 2n-1$. Therefore, this inequality can be regarded as a criterion of integrability. However, in this paper we explore another criterion of integrability given by  the so-called hydrodynamic reduction method.

\section{Gibbons--Tsarev systems}

Gibbons-Tsarev systems are the main ingredient of the approach to integrability of Whitham type hierarchies and, more generally, to integrability of quasi-linear systems of the form  (\ref{hydr}) based on hydrodynamic reductions. In this approach 
hydrodynamic reductions of a given hierarchy are parameterized by solutions of a Gibbons-Tsarev system. In this Section we explain a connection between  Gibbons-Tsarev systems and GT structures.

Let $p_1,..,p_M,~v_1,...,v_m$ be functions of auxiliary  variables $r_1,...,r_M$ and $\partial_i=\frac{\partial}{\partial r_i}$.

{\bf Definition 5.1.}
A Gibbons--Tsarev system  is a compatible system of partial differential equations of the form.
\begin{eqnarray}
\partial_ip_j=f(p_i,p_j,v_1,\dots,v_m)\partial_iv_1,\quad
i\ne j,\quad i,j=1,\dots,M, ~~~~~ \nonumber \\
\partial_iv_j=g_j(p_i,v_1,\dots,v_m)\partial_iv_1,\quad
j=2,\dots,m,\quad i=1,\dots,M,~~ \label{gtsys} \\
\partial_i\partial_jv_1=q(p_i,p_j,v_1,\dots,v_m)\partial_i v_1\partial_j v_1,\quad
i\ne j,\quad i,j=1,\dots,M.  \nonumber
\end{eqnarray}

{\bf Remark 5.1.} It follows from the  compatibility assumption that the space of solutions of a Gibbons-Tsarev system is  locally parameterized by $2M$ functions in one variable. Note that $f,~g_i,~ q$ do not depend on $M$ and therefore $M$ can be 
arbitrary large for a given Gibbons-Tsarev system.

We say that a Gibbons-Tsarev system is non-degenerate if  $f(p_1,p_2,v_1,...,v_m)$ has a pole of order one on the diagonal $p_2=p_1$. Assume in the sequel that all Gibbons-Tsarev systems are non-degenerate. 

{\bf Proposition 5.1.} There exists a one-to-one correspondence between non-degenerate  Gibbons--Tsarev systems  and local GT structures.

{\bf Proof.} Redefining $f,g_i$ from (\ref{gtsys}) we write a Gibbons-Tsarev system in the form
$$\partial_ip_j=\frac{f(p_i,p_j,v_1,\dots,v_m)}{g_1(p_i,v_1,...,v_m)}\partial_iv_1,\quad
i\ne j,\quad i,j=1,\dots,M,$$
\begin{equation}\label{gtsys1}
\frac{\partial_iv_1}{g_1(p_i,v_1,\dots,v_m)}=\frac{\partial_iv_j}{g_j(p_i,v_1,\dots,v_m)},\quad
j=2,\dots,m,\quad i=1,\dots,M,
\end{equation}
$$\partial_i\partial_jv_1=q(p_i,p_j,v_1,\dots,v_m)\partial_i v_1\partial_j v_1,\quad
i\ne j,\quad i,j=1,\dots,M$$
where $f(p_1,p_2)=\frac{1}{p_1-p_2}+O(1)$. Indeed, $\frac{1}{g_1(p_i)}$ is the residue of $f(p_i,p_j)$ from (\ref{gtsys}) at $p_j=p_i$. Write 
$$g(p)=\sum_{i=1}^mg_i(p,v_1,...,v_m)\frac{\partial}{\partial v_i}.$$
Compatibility of the system  (\ref{gtsys1})  implies  $\partial_1\partial_2\phi(p_3,v_1,...,v_m)=\partial_2\partial_1\phi(p_3,v_1,...,v_m)$ for an arbitrary function $\phi$. This can be written as
$$\Big(f(p_1,p_2)\frac{\partial}{\partial p_2}+f(p_1,p_3)\frac{\partial}{\partial p_3}+g(p_2)\Big)\Big((f(p_2,p_3)\frac{\partial}{\partial p_3}+g(p_2))\phi(p_3)\cdot\frac{\partial_2u_1}{g_1(p_2)}\Big)\cdot\frac{\partial_1u_1}{g_1(p_1)}=$$
$$\Big(f(p_2,p_1)\frac{\partial}{\partial p_1}+f(p_2,p_3)\frac{\partial}{\partial p_3}+g(p_1)\Big)\Big((f(p_1,p_3)\frac{\partial}{\partial p_3}+g(p_1))\phi(p_3)\cdot\frac{\partial_1u_1}{g_1(p_1)}\Big)\cdot\frac{\partial_2u_1}{g_1(p_2)}.~$$
Expanding this equation and equating coefficients at $\phi$ and $\phi_{p_3}$ we get
$$f(p_1,p_2)f(p_2,p_3)_{p_2}-f(p_2,p_1)f(p_1,p_3)_{p_1}+f(p_1,p_3)f(p_2,p_3)_{p_3}-f(p_2,p_3)f(f(p_1,p_3)_{p_3}+$$
$$+g(p_1)(f(p_2,p_3))-g(p_2)(f(p_1,p_3))+$$
$$f(p_2,p_3)\Big(g_1(p_1)\frac{\partial_1\partial_2u_1}{\partial_1u_1\partial_2u_1}-\frac{1}{g_1(p_2)}\Big(f(p_1,p_2)\frac{\partial}{\partial p_2}+g(p_1)\Big)(g_1(p_2))\Big)-$$
$$f(p_1,p_3)\Big(g_1(p_2)\frac{\partial_1\partial_2u_1}{\partial_1u_1\partial_2u_1}-\frac{1}{g_1(p_1)}\Big(f(p_2,p_1)\frac{\partial}{\partial p_1}+g(p_2)\Big)(g_1(p_1))\Big)=0,$$

$$f(p_1,p_2)g^{\prime}(p_2)-f(p_2,p_1)g^{\prime}(p_1)+[g(p_1),g(p_2)]+$$
$$\Big(g_1(p_1)\frac{\partial_1\partial_2u_1}{\partial_1u_1\partial_2u_1}-\frac{1}{g_1(p_2)}\Big(f(p_1,p_2)\frac{\partial}{\partial p_2}+g(p_1)\Big)(g_1(p_2))\Big)g(p_2)-$$
$$\Big(g_1(p_2)\frac{\partial_1\partial_2u_1}{\partial_1u_1\partial_2u_1}-\frac{1}{g_1(p_1)}\Big(f(p_2,p_1)\frac{\partial}{\partial p_1}+g(p_2)\Big)(g_1(p_1))\Big)g(p_1)=0.$$

Expanding the first of these equations near  the  diagonal $p_2=p_3$ and noting that
$$f(p_1,p_2)f(p_2,p_3)_{p_2}+f(p_1,p_3)f(p_2,p_3)_{p_3}-f(p_2,p_3)f(p_1,p_3)_{p_3}=-\frac{2f(p_1,p_2)_{p_2}}{p_2-p_3}+O(1)$$
we obtain
$$g_1(p_1)\frac{\partial_1\partial_2u_1}{\partial_1u_1\partial_2u_1}-\frac{1}{g_1(p_2)}\Big(f(p_1,p_2)\frac{\partial}{\partial p_2}+g(p_1)\Big)(g_1(p_2))=2f(p_1,p_2)_{p_2}.$$
Substituting this into our equations we arrive  at  relations  (\ref{GT2}),  (\ref{GT1})  for  a  local GT structure. 

One can check that all these steps are invertible and any local GT structure with relations  (\ref{GT1}),  (\ref{GT2}) gives a Gibbons-Tsarev system  (\ref{gtsys1}) with
$$\partial_1\partial_2u_1=\Big(\frac{2f(p_1,p_2)_{p_2}}{g_1(p_1)}+\frac{1}{g_1(p_1)g_1(p_2)}\Big(f(p_1,p_2)\frac{\partial}{\partial p_2}+g(p_1)\Big)(g_1(p_2))\Big)\partial_1u_1\partial_2u_1.$$   $\Box$

\section{Integrability of Whitham type hierarchies}

In this Section we explain a relation between integrable Whitham type hierarchies and enhanced  GT structures.

{\bf Proposition 6.1.} A  Whitham type hierarchy with potentials $h_i(p,v_1,...,v_m),~i=1,...,N$ is  integrable by hydrodynamic reductions  if and only if   there exists a Gibbons-Tsarev system (\ref{gtsys1}) such that 
\begin{equation}\label{int}
h_j^{\prime}(p_1)\partial_1(h_i(p_2))=h_i^{\prime}(p_1)\partial_1(h_j(p_2)),~~~i,j=1,...,N
\end{equation}
by virtue of (\ref{gtsys1}).

{\bf Proof.} The equation (\ref{int})  can be written as 
\begin{equation}\label{cr}
f(p_1,p_2)=\frac{\sum_{k=1}^m\Big(h_i^{\prime}(p_1)h_j(p_2)_{v_k}-
h_j^{\prime}(p_1)h_i(p_2)_{v_k}\Big)g_k(p_1)}
{h_j^{\prime}(p_1)h_i^{\prime}(p_2)-h_j^{\prime}(p_2)h_i^{\prime}(p_1)}
\end{equation}
and, therefore, coincides with the formula (77) from \cite{os11}.  It is proven in \cite{os11}  that the equation (\ref{cr})  is equivalent to  the  integrability of a given Whitham type hierarchy. $\Box$

{\bf Proposition 6.2.} There exists a one-to-one correspondence between integrable  Whitham type hierarchies and enhanced local GT structures. Under this correspondence the space of potentials of  a  Whitham type hierarchy  coincides with 
the space of solutions of the linear system  (\ref{GT5}).

{\bf Proof.} Write (\ref{int}) as $\frac{\partial_1(h_i(p_2))}{h_i^{\prime}(p_1)}=\frac{\partial_1(h_j(p_2))}{h_j^{\prime}(p_1)}$.  By executing  $\partial_1$ in numerators we get
$$
\frac{f(p_1,p_2)h_i^{\prime}(p_2)+g(p_1)(h_i(p_2))}{h_i^{\prime}(p_1)}=\frac{f(p_1,p_2)h_j^{\prime}(p_2)+g(p_1)(h_j(p_2))}{h_j^{\prime}(p_1)}.
$$
Let $\lambda(p_1,p_2)=\frac{f(p_1,p_2)h_i^{\prime}(p_2)+g(p_1)(h_i(p_2))}{h_i^{\prime}(p_1)}$, this function does not depend on $i$. Therefore, we get
$$g(p_1)(h_i(p_2))=\lambda(p_1,p_2)h_i^{\prime}(p_1)-f(p_1,p_2)h_i^{\prime}(p_2)$$
which coincides with  (\ref{GT5}). Applying the relation  (\ref{GT1})  to $h_i(p_3)$ we can write
$$g(p_1)g(p_2)(h_i(p_3))-g(p_2)g(p_1)(h_i(p_3))=$$$$f(p_2,p_1)b^{\prime}(p_1)(h_i(p_3))-f(p_1,p_2)b^{\prime}(p_2)(h_i(p_3))+2f(p_2,p_1)_{p_1}b(p_1)(h_i(p_3))-2f(p_1,p_2)_{p_2}b(p_2)(h_i(p_3)).$$
Computing the l.h.s. and the r.h.s. of this relation by virtue of (\ref{GT5})  we obtain  (\ref{GT4}).  $\Box$

\section{The universal  Whitham  hierarchy}

In this Section we use notations introduced in Section 2,  including $G(p)$ and  $F(p_1,p_2)$.

According to \cite{kr2} the universal  Whitham  hierarchy  is given by potentials obtained by integration of meromorphic differentials on a Riemann surface. We are going to construct such  an  hierarchy explicitly\footnote{We need to choose  
constants of integrations carefully  in order to obtain an integrable hierarchy.} and prove that it is integrable by 
hydrodynamic reductions. 

{\bf Proposition 7.1.} Fix constants $s_1,...,s_m$ such that $s_1+...+s_m=1$ (the simplest  possibility is $m=1$ and $s_1=1$).  The following formulas define an enhanced GT structure:
 \begin{equation}\label{uwh}
g(p)=\sum_{j=1}^nF(p,u_j)\frac{\partial}{\partial u_j}+\sum_{k=1}^mF(p,w_j)\frac{\partial}{\partial w_k}+G(p),~~~f(p_1,p_2)=F(p_1,p_2),
\end{equation}
$$\lambda(p_1,p_2)=\frac{E(p_1,p_2)_{p_1}}{E(p_1,p_2)}-\sum_{k=1}^ms_k\frac{E(p_1,w_k)_{p_1}}{E(p_1,w_k)}.$$
Moreover, the following functions belong to the space of potentials  of this enhanced GT structure:
$$h_j(p)-h_1(p),~j=2,...,n,~~~q_{\alpha}(p)-\sum_{k=1}^ms_kq_{\alpha}(w_k),~\alpha=1,...,g$$
where
\begin{equation}\label{pot1}
h_j(p)=\ln(E(p,u_j))-\sum_{k=1}^ms_k\ln(E(u_j,w_k)).
\end{equation}

{\bf Proof.} We need to prove identities (\ref{GT4}) and (\ref{GT5}) for given $\lambda(p_1,p_2)$ and potentials. This can be done by straightforward computation using identities from Proposition 2.2. The simplest way is to start  from identity 
(\ref{GT5}) for $h_j(p)-h_1(p)$ and check it using identity (\ref{GTc4}). It is clear from the proof of the Proposition 6.2 that (\ref{GT4}) is a consequence of (\ref{GT5}). It follows from Proposition 3.3 that $\frac{1}{2\pi i}\int_{b_{\alpha}}\lambda(t,p)d t$ are 
also potentials of our hierarchy. Computing these integrals by virtue of (\ref{aut})  we conclude that the functions  $q_{\alpha}(p)-\sum_{k=1}^ms_kq_{\alpha}(w_k),~\alpha=1,...,g$ belong to the space of potentials.  $\Box$

{\bf Proposition 7.2.} The universal Whitham hierarchy is integrable by hydrodynamic reductions.

{\bf Proof.} It is clear that the vector space spanned  by  derivatives with pespect to $p$  of potentials described in the previous Proposition coincides with the space of meromorphic differentials on  $\mathcal{E}$ holomorphic outside $u_1,...,u_n$ and with poles of 
order  less or equal to one   in these points. Therefore, we obtain a part of the universal Whitham hierarchy. In order to obtain the full hierarchy we apply the procedure of colliding point, see Proposition 3.2 and Remark 3.4. This proves the  Proposition 
in the case $g>1$.

In the case $g=0$ we define an enhanced GT structure by (\ref{g0}) and set $\lambda(p_1,p_2)=\frac{1}{p_1-p_2}$. The space of potentials contains the functions  $h_j(p)-h_1(p),~j=2,...,n+2$ where $h_j(p)=\ln(p-u_j),~j=1,...,n$, $h_{n+1}(p)=\ln(p)$ and $h_{n+2}(p)=\ln(p-1)$. This gives a part of the universal Whitham hierarchy corresponding to meromorphic differentials on $\C P^1$ with poles of order   less or equal to one  in $u_1,...,u_n,0,1$. To obtain the  full hierarchy we collide these  points 
by a procedure similar to the  one in the proof  of Proposition 3.2, see also Remarks 3.3 and  3.4.

In the case $g=1$ we define an enhanced GT structure by (\ref{g1}) and set $$\lambda(p_1,p_2)=\rho(p_1-p_2,\tau)-\rho(p_1)-2\pi i.$$ The space of potentials contains $p-\tau$ and  the functions $h_j(p)-h_1(p),~j=2,...,n$ 
where $h_j(p)=\ln(\theta(p-u_j,\tau))-\ln(\theta(u_j,\tau))$. This gives a part of the universal Whitham hierarchy corresponding to meromorphic differentials on $\mathcal{E}$ with poles of order less or equal to one in $u_1,...,u_n$. 
To obtain the full hierarchy we collide some of  these  points by a procedure similar to one in the proof  of Proposition 3.2, see also Remark 3.4.  $\Box$

\bigskip

\vskip.4cm
\noindent
{\bf Acknowledgments.}

I am grateful to Maxim Kontsevich for useful discussions. This paper was completed during my visit to IHES. I am grateful to IHES for hospitality and an  excellent working atmosphere.

\end{document}